\documentclass[12pt]{amsart}
\usepackage{amssymb,latexsym, amscd, a4wide}
\usepackage[all]{xy}
\usepackage{graphicx}
\usepackage{mathrsfs}
\usepackage{amsmath}
\usepackage{pb-diagram}
\usepackage{bbold}
\usepackage{color}
\usepackage{amssymb}

\newtheorem{thm}{Theorem}[section]

\newtheorem{prop}[thm]{Proposition}
\newtheorem{cor}[thm]{Corollary}
\newtheorem{ex}[thm]{Example}

\newtheorem{rem}[thm]{Remark}

\newtheorem*{thm*}{Theorem}
\newtheorem*{prop*}{Proposition}

 \numberwithin{equation}{section}
 \newcommand{\nc}{\newcommand}
\nc{\ssn}{\subsection{}} \nc{\sssn}{\subsubsection{}}

\newcommand{\dem }{\noindent\textbf{Proof}. }
\newcommand{\findem }{\hfill $\Box$ \vskip0.5cm }

\numberwithin{equation}{section}
\newcommand{\Ind}{\operatorname{Ind}}

\nc{\oR}{\ol{R}}
\newcommand{\ol}{\ensuremath{\overline}}

\nc{\id}{\textrm{id}}

\newcommand{\Ext}{\operatorname{Ext}}

\newcommand{\Hom}{\operatorname{Hom}}

\newcommand{\ad}{\hbox{\ensuremath{\operatorname{ad}}}}

\newcommand{\Z}{\ensuremath{\mathbb{Z}}}

\newcommand{\C}{\ensuremath{\mathbb{C}}}

\newcommand{\g}{\ensuremath{\mathfrak{g}}}
\newcommand{\h}{\ensuremath{\mathfrak{h}}}
\renewcommand{\b}{\ensuremath{\mathfrak{b}}}

\newcommand{\n}{\ensuremath{\mathfrak{n}}}

\newcommand{\mc}[1]{\mathcal{#1}} 
\newcommand{\mf}[1]{\mathfrak{#1}} 

\begin{document}
\author{Juan Camilo Arias, Vyacheslav Futorny and Andr\'e de Oliveira}
\title{The category of reduced imaginary Verma modules}
\address{Institute of Mathematics and Statistics, University of S\~ao Paulo, S\~ao Paulo, BRAZIL. }
\email{jcarias@ime.usp.br }
\address{Shenzhen International Center for Mathematics, Southern University of Science and Technology, China and Institute of Mathematics and Statistics, University of S\~ao Paulo, S\~ao Paulo, BRAZIL.}
\email{vfutorny@gmail.com}
\address{Institute of Mathematics and Statistics, University of S\~ao Paulo, S\~ao Paulo, BRAZIL. }
\email{andre2.oliveira@usp.br}

\keywords{}

\maketitle

\begin{abstract} For an arbitrary affine Lie algebra  we study an analog of the category $\mathcal O$ for the natural Borel subalgebra and zero central charge. We show that such category 
is semisimple having
 the reduced imaginary Verma modules as its simple objects. This generalizes  the result of Cox, Futorny, Misra in the case of affine $\mf{sl}_2$.
\end{abstract}

\section{Introduction}

Let $A=(a_{ij})_{0\leq i,j\leq N}$ be a generalized affine Cartan matrix over $\C$ with associated affine Lie algebra $\hat{\g}$ and Cartan subalgebra $\hat{\h}$. 
Let $\Pi =\{\alpha_0, \alpha_1, \cdots, \alpha_N\}$ be the set of simple roots, $\delta$ the indivisible imaginary root and $\Delta$ the root system of $\hat{\g}$. A subset $S\subseteq \Delta$ is a closed partition if for any $\alpha, \beta \in S$ and $\alpha + \beta \in \Delta$ then $\alpha + \beta \in S$, $\Delta = S \cup (-S) $ and $S\cap (-S) = \emptyset$.  The classification of closed partitions for root system of affine Lie algebras was obtained by H. Jakobsen and V. Kac in \cite{JK01} and \cite{JK02} and independently by V. Futorny in \cite{F03} and \cite{F04}. They show that closed partitions are parameterized by subsets $X\subseteq\Pi$ and that (contrary to what happens in the finite case) there exists a finite number (greater than 1) of inequivalent Weyl group orbits of closed partitions. When $X=\Pi$ we get that $S=\Delta_{+}$ and we can developed the standard theory of Verma modules, but in the case $X\subsetneq \Pi$ we obtain new Verma-type modules called non-standard Verma modules. \\

The theory of non-standard Verma modules was initiated by V. Futorny in \cite{F05} (see also \cite{F01}) in the case $X=\varnothing$ and continued by B. Cox in \cite{C01} for arbitrary $X\subsetneq \Pi$. The case $X=\varnothing$ give rise to the natural Borel subalgebra associated to the natural partition $\Delta_{\textrm{nat}} = \{ \alpha + n\delta \ | \ \alpha\in \Delta_{0,+} \ , \ n  \in \Z   \} \cup \{ k\delta \ | \ k \in \mathbb{Z}_{>0}   \}$. The Verma module $M(\lambda)$, of highest weight $\lambda$, induced by the natural Borel subalgebra is called imaginary Verma module for $\hat{\g}$, when it is not irreducible it has an irreducible quotient called reduced imaginary Verma module. Unlike the standard Verma modules, imaginary Verma modules contain both finite and infinite dimensional weight spaces. Similar results hold for more general non-standard Verma modules.\\

In \cite{CFM03}, while studying crystal bases for reduced imaginary Verma modules of $\hat{\mf{sl}_2}$, it was consider a suitable category of modules, denoted $\mc{O}_{red,im}$, with the properties that any module in this category is a reduced imaginary Verma module or it is a direct sum of these modules. In this paper, by appropriate modifications we first define a category $\mc{O}_{red,im}$ for any affine Lie algebra and we show that all irreducible modules in this category are reduced imaginary Verma modules and, moreover, that any arbitrary module in $\mc{O}_{red,im}$ is a direct sum of reduced imaginary Verma modules. \\

{\it It should be noted that the results presented in this paper hold for both untwisted and twisted affine Lie algebras}.
The paper is organized as follows. In Sections $2$ and $3$, we define, set the notations and summarize the basic results for affine algebras, closed partitions and imaginary Verma modules. In section $4$ we introduce the category $\mc{O}_{red,im}$ and present some of its properties. Finally, in section $5$ we present the main results of this paper. 

\section{Preliminaries}

In this section we fixed some notation and the preliminaries about affine algebras and root datum are set up. 

\subsection{Affine algebras} 

Let $A=(a_{ij})_{0\leq i,j\leq N}$ be a generalized affine Cartan matrix over $\C$ with associated affine Lie algebra $\hat{\g}$. Let $D=diag(d_0, \ldots, d_N)$ be a diagonal matrix with relatively primes integer entries such that $DA$ is symmetric. The Lie algebra $\hat{\g}$ has a Chevalley-Serre presentation given by generators $e_i, f_i, h_i$ for $0\leq i \leq N$ and $d$ which are subject to the defining relations:

$$ [h_i,h_j]=0 \quad [d,h_i]=0 \quad [h_i,e_j]=a_{ij}e_j \quad [h_i,f_j]=-a_{ij}f_j $$
$$[e_i,f_j]=\delta_{i,j}h_i \quad [d,e_i]=\delta_{0,i}e_i \quad [d,f_i]=-\delta_{0,i}f_i$$
$$(\ad e_i)^{1-a_{ij}}(e_j)=0 \quad (\ad f_i)^{1-a_{ij}}(f_j)=0$$

Let $\hat{\h}$ be the Cartan subalgebra of $\hat{\g}$ which is the span of $\{h_0, \ldots, h_N,d\}$. \\

Recall that affine Lie algebras are classified into two classes: untwisted and twisted, 
see \cite[Ch. 6-8]{Kac}. In the untwisted case, $\hat{\g}$ has a natural realization known as loop space realization which is defined by 

$$\hat{\g} = \g\otimes \C[t,t^{-1}]\oplus \C c \oplus \C d$$
\vspace{3mm}

where $\g$ is the simple finite dimensional Lie algebra with Cartan matrix $(a_{ij})_{1\leq i,j \leq N}$, $c$ is a central element, $d$ is a degree derivation such that $[d,x\otimes t^n]=nx\otimes t^n$ for any $x\in \g$ and $n\in \Z$ and we have $[x\otimes t^n, y\otimes t^m] = [x,y]\otimes t^{n+m} + \delta_{n,-m}n(x|y)c$ for all $x,y\in\g$, $n,m\in\Z$ where $(-|-)$ is a symmetric invariant bilinear form on $\g$.\\

On the other hand, twisted affine Lie algebras are described as fixed points of automorphisms of untwisted algebras. Concretely, let $\tilde\mu$ be an automorphism of order $r=2$ or $r=3$ of the Coxeter-Dynkin diagram of $\g$ and let $\overline\mu$ be the corresponding diagram automorphism of $\g$.\\

Then $\overline\mu$ can be extended to an automorphism $\mu$ on $\hat{\g} = \g\otimes \C[t,t^{-1}]\oplus \C c \oplus \C d$ defined as $\mu(x \otimes t^{m}) = (-1)^{m} (\overline{\mu}(x) \otimes t^{m})$, for $x \in \g$, $m \in \mathbb{Z}$, $\mu(c) = c$, $\mu(d) = d$ and extended by linearity. The \textit{twisted affine Lie algebra} $(\hat{\g})^{\mu}$ is the subalgebra of fixed points of $\mu$.\\

For example, when $r = 2$,
$$(\hat{\g})^{\mu} = \left(\sum_{m \in \mathbb{Z}}\mu_{0} \otimes t^{2m}\right) \oplus \left(\sum_{m \in \mathbb{Z}}\mu_{1} \otimes t^{2m+1}\right) \oplus \mathbb{C}c \otimes \mathbb{C}d$$
where $\mu_{0} = \{x \in \g \ | \ \overline{\mu}(x) = x\}$ and $\mu_{1} = \{x \in \g \ | \ \overline{\mu}(x) = -x\}$ (see \cite{F02}).\\

\subsection{Root datum and closed partitions}

Let $I_{0} = \{1, \ldots, N\}$ and $\Delta_0$ be the root system of $\g$ with $\theta$ being the longest positive root. We denote by $Q_0$ and $P_0$ the root and weight lattices of $\g$. Let $I=\{0,1,\ldots, N\}$, $\Delta$ the root system of $\hat{\g}$ with simple roots $\Pi=\{\alpha_0, \alpha_1, \ldots, \alpha_N\}$ and let $\delta=\alpha_0+\theta$ be the  \textit{indivisible imaginary root}. $Q$ denotes the root lattice, $P$ the weight lattice, and $\check{Q}$, $\check{P}$ denotes the coroot and coweight lattices, respectively. $\Delta^\textrm{re}$ and $\Delta^\textrm{im}$ denotes the real and the imaginary sets of roots for $\Delta$.\\

A subset $S$ of $\Delta$ is said to be closed if whenever $\alpha, \beta \in S$ and $\alpha + \beta \in \Delta$ then $\alpha + \beta \in S$. We also say that $S$ is a closed partition if $S$ is closed, $\Delta = S \cup (-S) $ and $S\cap (-S) = \emptyset$. Closed partitions were classified in \cite{F03} and \cite{F04} (see also \cite{JK01} and \cite{JK02}).\\

For an untwisted affine Lie algebra $\hat{\g}$, there are two interesting closed partitions of the root system $\Delta$, the \textit{standard partition} and the \textit{natural partition}, which give rise to two distinct Borel subalgebras that are not conjugate. \\

The standard partition is defined by 

$$ \Delta_{\textrm{st}} = \{ \alpha + n\delta \ | \ \alpha\in \Delta_{0} \ , \ n  \in \Z_{>0}   \} \cup \Delta_{0,+} \cup \{ k\delta \ | \ k \in \mathbb{Z}_{>0}   \}$$

 and the natural partition by 
 
 $$ \Delta_{\textrm{nat}} = \{ \alpha + n\delta \ | \ \alpha\in \Delta_{0,+} \ , \ n  \in \Z   \} \cup \{ k\delta \ | \ k \in \mathbb{Z}_{>0}   \}$$

 The respective Borel subalgebras, called \textit{standard Borel subalgebra} and \textit{natural Borel subalgebra}, are defined by
 $$\b_{\textrm{st}} = \Big( \g \otimes t\mathbb{C}[t]\Big) \oplus \n \oplus \h \oplus \mathbb{C}c \oplus \mathbb{C}d$$ 
 and
 $$\b_{\textrm{nat}} = \Big( \n \otimes \mathbb{C}[t,t^{-1}]\Big) \oplus \Big( \h \otimes t\mathbb{C}[t]\Big) \oplus \h \oplus \mathbb{C}c \oplus \mathbb{C}d$$ 

where $\mf{n} = \bigoplus_{\alpha\in \Delta_{0,+}} \g_{\alpha}$, is the nilpotent Lie subalgebra of the finite Lie algebra $\g$. \\
 
As already mentioned above, a twisted affine algebra is a fixed point set in $\hat{\g}$ of a non-trivial symmetry of Chevalley generators and, in this case,  $\b_{\textrm{nat}}$ is the intersection of the fixed point set with the natural Borel subalgebra of $\hat{\g}$. For more details see \cite{JK01}.\\

In this paper, we are going to work with the natural partition of the root system $\Delta_{\textrm{nat}}$. \\

\section{Imaginary Verma modules}

Let $S$ be a closed partition of the root system $\Delta$. Let $\hat{\g}$ be the untwisted affine Lie algebra which has, with respect to the partition $S$, the triangular decomposition $\hat{\g}=\hat{\g}_S\oplus \hat{\h}\oplus \hat{\g}_{-S}$, where $\hat{\g}_S = \bigoplus_{\alpha\in S}\hat{\g}_{\alpha}$ and $\hat{\h} = \h \oplus \mathbb{C}c\oplus \mathbb{C}d$ is an affine Cartan subalgebra. Let $U(\hat{\g}_S)$ and $U(\hat{\g}_{-S})$ be, respectively, the universal enveloping algebras of $\hat{\g}_S$ and $\hat{\g}_{-S}$.\\

Let $\lambda\in P$. A weight $U(\hat{\g})$-module V is called an $S$-highest weight module with highest weight $\lambda$ if there is some non-zero vector $v\in V$ such that:
\vspace{3mm}
\begin{itemize}
    \item $u\cdot v = 0$ for all $u \in \hat{\g}_S$.
    \item $h \cdot v = \lambda(h)v$ for all $h \in \hat{\h}$.
    \item $V=U(\hat{\g})\cdot v \cong U(\hat{\g}_{-S})\cdot v$.
\end{itemize}
\vspace{3mm}

In what follows, let us consider $S$ to be the natural closed partition of $\Delta$, i.e., $S=\Delta_{\textrm{nat}}$ and so $\mf{b}_{\textrm{nat}} = \hat{\g}_{\Delta_{\textrm{nat}}} \oplus \hat{\mf{h}}$. We make $\C$  into a 1-dimensional $U(\mf{b}_{\textrm{nat}})$-module by picking a generating vector $v$ and setting $(x+h)\cdot v = \lambda(h)v$, for all $x\in \hat{\g}_{\Delta_{\textrm{nat}}}$ and $h\in \hat{\h}$. The induced module 

$$ M(\lambda) = U(\hat{\g})\otimes_{U(\mf{b}_{\textrm{nat}})}\C v \cong U(\hat{\g}_{-\Delta_{\textrm{nat}}})\otimes \C v  $$

\vspace{2mm}

is called an {\it imaginary Verma module} with $\Delta_{\textrm{nat}}$-highest weight $\lambda$. Equivalently, we can define $M(\lambda)$ as follows: Let $I_{\Delta_{\textrm{nat}}}(\lambda)$ the ideal of $U(\hat{\g})$ generated by $e_{ik}:= e_i\otimes t^k$, $h_{il}:= h_i\otimes t^l$ for $i\in I_0$, $k\in \Z$, $l \in \mathbb{Z}_{>0}$, and by $h_i - \lambda(h_i)\cdot 1$, $d-\lambda(d)\cdot 1$ and $c-\lambda(c)\cdot 1$. Then $M(\lambda) = U(\hat{\g})/I_{\Delta_{\textrm{nat}}}(\lambda)$.\\

The main properties of this modules, which hold for any affine Lie algebra, were proved in \cite{F01} (see also \cite{F06} for more properties on this modules), we summarize them in the following. 

\begin{prop} Let $\lambda \in P$ and let $M(\lambda)$ be the imaginary Verma module of $\Delta_{\textrm{nat}}$-highest weight $\lambda$. Then $M(\lambda)$ has the following properties:
\begin{enumerate}
    \item The module $M(\lambda)$ is a free $U(\hat{\g}_{-\Delta_{\textrm{nat}}})$-module of rank 1 generated by the $\Delta_{\textrm{nat}}$-highest weight vector $1\otimes 1$ of weight $\lambda$.
    \item $M(\lambda)$ has a unique maximal submodule.
    \item Let $V$ be a $U(\hat{\g})$-module generated by some $\Delta_{\textrm{nat}}$-highest weight vector $v$ of weight $\lambda$. Then there exists a unique surjective homomorphism $\phi: M(\lambda) \to V$ such that $1\otimes 1 \mapsto v$. 
    \item $\dim M(\lambda)_{\lambda} = 1$. For any $\mu=\lambda-k\delta$, $k\in \Z_{>0}$, $0<\dim M(\lambda)_\mu < \infty$. If $\mu \neq \lambda - k\delta$ for any integer $k\geq 0$ and $M(\lambda)_\mu\neq 0$, then $\dim M(\lambda)_\mu = \infty$.
    \item Let $\lambda, \mu \in \hat{\h}^*$. Any non-zero element of $\Hom_{U(\hat{\g})}(M(\lambda), M(\mu))$ is injective. 
    \item The module $M(\lambda)$ is irreducible if and only if $\lambda(c)\neq 0$.
\end{enumerate}
\end{prop}

Suppose now that $\lambda(c)=0$ and consider the ideal $J_{\Delta_{\textrm{nat}}}(\lambda)$ generated by $I_{\Delta_{\textrm{nat}}}(\lambda)$ and $h_{il}$, $i\in I_0$ and $l\in \Z\setminus\{0\}$. Set 

$$\tilde{M}(\lambda) = U(\hat{\g})/J_{\Delta_{\textrm{nat}}}(\lambda)$$

Then $\tilde{M}(\lambda)$ is a homomorphic image of $M(\lambda)$ which we call {\it reduced imaginary Verma module}. The following is proved in \cite{F01}, Theorem 1.

\begin{prop}
$\tilde{M}(\lambda)$ is irreducible if and only if $\lambda(h_i)\neq 0$ for all $i\in I_0$.
\end{prop}

\section{The category $\mc{O}_{red,im}$}

Consider the Heisenberg subalgebra $G$ which by definition is 
$$G= \bigoplus_{k\in \Z\setminus\{0\}} \hat{\g}_{k\delta} \oplus \C c$$

We will say that a $\hat{\g}$-module V is $G$-compatible if:

\begin{enumerate}
    \item [(i)] $V$ has a decomposition $V=T(V)\oplus TF(V)$ where $T(V)$ and $TF(V)$ are non-zero $G$-modules, called, respectively, torsion and torsion free module associated to $V$.
    \item [(ii)] $h_{im}$ for $i\in I_0$, $m\in \Z\setminus\{0\}$ acts bijectively on $TF(V)$, i.e., they are bijections on $TF(V)$.
    \item [(iii)] $TF(V)$ has no non-zero $\hat{\g}$-submodules.
    \item [(iv)] $G\cdot T(V)=0$.
\end{enumerate}

Consider the set 

$$ \hat{\h}^*_{red} = \{ \lambda \in \hat{\h}^* \; | \;  \lambda(c)=0, \lambda(h_i)\notin\Z_{\geq 0} \mbox{ for any } i\in I_0 \} $$

We define the category $\mc{O}_{red,im}$ as the category whose objects are $\hat{\g}$-modules $M$ such that 

\begin{enumerate}
    \item $M$ is $\hat{\h}^*_{red}$-diagonalizable, that means,  $$ M = \bigoplus_{\nu \in \hat{\h}^*_{red}} M_{\nu}, \mbox{  where  } M_{\nu} = \{ m\in M | h_im=\nu(h_i)m, dm = \nu(d)m, i\in I_0 \}$$
    \item For any $i\in I_0$ and any $n\in \Z$, $e_{in}$ acts locally nilpotently.
    \item $M$ is $G$-compatible.
    \item The morphisms between modules are $\hat{\g}$-homomorphisms
\end{enumerate}

\begin{ex}
    Reduced imaginary Verma modules belongs to $\mc{O}_{red,im}$. Indeed, for $\tilde{M}(\lambda)$ consider $T(\tilde{M}(\lambda)) = \C v_{\lambda}$ and $TF(V) = \bigoplus_{k\in\Z, n_1, \ldots, n_N\in\Z_{\geq0}} \tilde{M}(\lambda)_{\lambda+k\delta - n_1\alpha_1 - \ldots -n_N\alpha_N}$, and at least one $n_j\neq 0$. Moreover, direct sums of reduced imaginary Verma modules belongs to $\mc{O}_{red,im}$.
\end{ex}

Recall that a loop module for $\hat{\g}$ is any representation of the form $\hat{M} := M \otimes \mathbb{C}[t,t^{-1}]$ where $M$ is a $\mathfrak{g}$-module and the action of $\hat{\g}$ on $\hat{M}$ is given by 
$$(x \otimes t^{k})(m \otimes t^{l}) := (x \cdot m) \otimes t^{k+l} \quad , \quad c(m \otimes t^{l}) = 0$$
for $x \in \mathfrak{g}$, $m \in M$ and $k,l \in \mathbb{Z}$. Here $x \cdot m$ is the action of $x \in \mathfrak{g}$ on $m \in M$.

\begin{prop}
  Let $M$ is a $\mathfrak{g}$-module in the BGG category $\mathcal{O}$.  Then the  loop module $\hat{M}$ can not lie in $\mc{O}_{red,im}$.
\end{prop}

\dem Let $M \in \mathcal{O}$ and let $\hat{M}$ be its associated loop module. If $M$ is finite dimensional, it is a direct sum of finite dimensional irreducible $\g$-modules, and these have highest weights which are non-negative integers when evaluated in $h_i$ for any $i\in I_0$. So, condition (1) is not satisfied and $\hat{M}$ does not belongs to $\mc{O}_{red,im}$. 
 Assume now that $M$ is an infinite dimensional $\g$-module. Note that condition (2) is satisfied as 
 $\n$ acts locally nilpotently on $M$. 
 If condition (1) does not hold, we are done. 
 Suppose that (1) holds and that $\hat{M}$ is $G$-compatible. We have $\hat{M} = T(\hat{M}) \oplus TF(\hat{M})$ satisfying (i) - (iv) above. Take any nonzero element $\sum_{i=-k}^{k}m_{i} \otimes t^{i} \in T(\hat{M})$ with $m_{i} \in M_{\mu}$ for some weight $\bar{\mu} \in \hat{\h}^*_{red}$. Then by (iv) we have
$$0 = (h_{j} \otimes t^{r})\left(\sum_{i=-k}^{k}m_{i} \otimes t^{i} \right) = \sum_{i=-k}^{k}(h_{j} \cdot m_{i}) \otimes t^{i+r} = \bar{\mu}(h_{j})\left(\sum_{i=-k}^{k}m_{i} \otimes t^{i+r} \right)$$
where $ j \in I_0$, $r \in \mathbb{Z}\setminus\{0\}$. Hence $\bar{\mu}(h_{j}) = 0$, for any $ j \in I_0$, which contradicts to the fact that $\bar{\mu} \in \hat{\h}^*_{red}$. Then $T(\hat{M}) = 0$ and $\hat{M} = TF(\hat{M})$ which is a $\hat{\g}$-module contradicting (i) and (iii), and thus (3). This completes the proof.

\findem

\section{Main results}

In this section we will show that the category $\mc{O}_{red,im}$ is a semisimple category having reduced imaginary Verma modules as its simple objects. First we will show that reduced imaginary Verma modules have no nontrivial extensions in $\mc{O}_{red,im}$.

\begin{thm}
    If $\lambda,\mu \in \hat{\h}^*_{red}$ then $\Ext_{\mc{O}_{red,im}}^1(\tilde{M}(\lambda), \tilde{M}(\mu)) = 0$.
\end{thm}

\dem Let $M$ be an extension of $\tilde{M}(\lambda)$ and $\tilde{M}(\mu)$ that fits in the following short exact sequence

\[ \xymatrix{0 \ar[r] & \tilde{M}(\lambda) \ar[r]^{\iota} & M \ar[r]^{\pi} & \tilde{M}(\mu) \ar[r] & 0 } \]

Suppose $\mu = \lambda +k\delta- \sum_{i=1}^N s_i\alpha_i $, for $s_i\in\Z$ and $k\in \Z$, and all $s_i$'s have the same sign or equal to $0$. First, consider the case when $s_i=0$ for all $i\in I_0$. Then $\mu = \lambda + k\delta$ and so, in $M$ there will be two  vectors  $v_{\lambda}$ and $v_{\mu}$ of weights $\lambda$ and $\mu$ respectively, annihilated by $\n\otimes \C[t, t^{-1}]$. Moreover, because of the condition (iv) in the definition of $G$-compatibility, these two points are isolated. So, $v_{\lambda}$ and $v_{\mu}$ are highest weight vectors, each of which generates an irreducible subrepresentation  
(isomorphic to $\tilde{M}(\lambda)$ and $\tilde{M}(\mu)$ respectively), and
the extension splits. Hence, we can assume that not all $s_i$ are equal to zero and that the map $\iota: \tilde{M}(\lambda) \to M$  in the short exact sequence is an inclusion. Assume that $s_i\in\Z_{\geq 0}$ for all $i$.\\

Let $\overline{v}_{\mu}\in M$ be a preimage under the map $\pi$ of a highest weight vector $v_{\mu}\in \tilde{M}(\mu)$ of weight $\mu$. We have $(\n\otimes \C[t, t^{-1}])v_{\mu}=Gv_{\mu}=0$, and 
 we are going to show that $G\overline{v}_{\mu}=0$.  Assume that $\overline{v}_{\mu}\notin T(M)$. Then we claim that $T(M)=\C v_{\lambda}$. Indeed, we have $\C v_{\lambda}\subset T(M)$. If $ u\in T(M)\setminus \C v_{\lambda}$ is some nonzero weight element, then $G\cdot u=0$ and $\pi(u)$ belongs to $T(\tilde{M}(\mu))=\C v_{\mu}$. If $\pi(u)=0$ then $u\in \tilde{M}(\lambda)$ which is a contradiction. If $\pi(u)$ is a nonzero multiple of $v_{\mu}$, then $u$ has weight $\mu$ and thus $u$ is a multiple of $\overline{v}_{\mu}$ which is again a contradiction. So, we assume $T(M)=\C v_{\lambda}$.\\

 Note that for any $i\in I_0$ and $m\in \Z\setminus \{0\}$ we have 
 $$\pi (h_{im}  \overline{v}_{\mu})=h_{im}  \pi (\overline{v}_{\mu})= h_{im}  v_{\mu} = 0.$$ Then $h_{im}  \overline{v}_{\mu}\in \tilde{M}(\lambda)$. Suppose there exists $j\in I_0$ such that $h_{jm}  \overline{v}_{\mu} \neq 0$ for $m\in \Z\setminus \{0\}$.  Because $h_{jm}  \overline{v}_{\mu} \in \tilde{M}(\lambda)$ and has weight $\mu+m\delta$, it belongs to $TF(\tilde{M}(\lambda))$. Hence, there exists a nonzero
  $v'\in \tilde{M}(\lambda)$ of weight $\mu$ such that $h_{jm} \overline{v}_{\mu} = h_{jm}  v'$. 
  Hence, $h_{jm} (\overline{v}_{\mu} - v')=0$ implying $\overline{v}_{\mu} - v' \in T(M) \cong \C v_{\lambda}$. Then $\overline{v}_{\mu} - v' = p \ v_{\lambda}$, for some $p \in \mathbb{C}$. 
Comparing the weight we arrive to a contradiction. Hence,  $h_{in} \overline{v}_{\mu}=0$. So, we get $G\overline{v}_{\mu}=0$.\\

%

Recall that the operators $e_{im}$ acts locally nilpotently on $\tilde{M}(\lambda)$.  We claim that 
$e_{im}\overline{v}_{\mu}=0$ for all possible $i$ and $n$. Indeed, assume that $e_{jm}\overline{v}_{\mu}\neq 0$ for some 
$j\in I_0$ and some integer $m$. Then $e_{im}\overline{v}_{\mu}\in \tilde{M}(\lambda)$. Consider the $\hat{\mf{sl}}_2$-subalgebra $\mf{s}(j)$ generated by $f_{jn}, e_{jn}$ and $h_{jl}$ for $n,l\in \Z$.
 Let $M_j$ be 
 an $\mf{s}(j)$-submodule of $M$ generated by $\overline{v}_{\mu}$. Then $M_j$ is an extension of 
  reduced imaginary Verma $\mf{s}(j)$-modules, one of which of highest weight $\mu$. Since $M\in \mc{O}_{red,im}$, we immediately see that $M_j$ is an object of the corresponding reduced category 
$\mc{O}_{red,im}(\mf{s}(j))$ for $\mf{s}(j)$. But this category is semisimple by \cite{CFM03}. Hence,
$e_{im}\overline{v}_{\mu}=0$ for all $i$ and $m$. Therefore, $\overline{v}_{\mu}$ generates a $\mathfrak{g}$-submodule of $M$ isomorphic to $\tilde{M}(\mu)$ and the short exact sequence splits. \\

 Assume now that $s_i\in\Z_{\leq 0}$ for all $i$ and not all of them are $0$. As $\tilde{M}(\mu)$ is irreducible and $\tilde{M}(\lambda)$ is a $\mathfrak{g}$-submodule of $M$, the short exact sequence splits completing the proof.

\findem

\begin{rem}
    Observe that modules $\tilde{M}(\lambda)$ and $\tilde{M}(\lambda-k\delta)$ have a nontrivial extension in the category of $\hat{\g}$-modules for any integer $k$. 
\end{rem}

\begin{thm}

    If $M$ is an irreducible module in the category $\mc{O}_{red,im}$, then $M\cong \tilde{M}(\lambda)$ for some $\lambda \in \hat{\mf{h}}_{red}^*$.
\end{thm}

\dem Let $M$ be an irreducible module in $\mc{O}_{red,im}$. As a $G$-module, $M\cong T(M)\oplus TF(M)$ where both summands are non-zero. Let $v\in T(M)$ be a non-zero element of weigh $\lambda \in \hat{\mf{h}}_{red}^*$. Then $h_{im}v=0$ for all $i\in I_0$ and all $m\in \Z\setminus \{0\}$. For each $i\in I_0$ let $p_i \in \Z_{>0}$ be the minimum possible integer such that $e_{i0}^{p_i}v=0$. If all $p_i=1$ we have $e_{i0}v=0$ and then, because $[h_{in},e_{i0}]=2e_{in}$ we get that $e_{in}v=0$ for all $i\in I_0$ and $n\in\Z\setminus \{0\}$. Hence, we have an epimorphism $\tilde{M}(\lambda) \twoheadrightarrow M$, since $\lambda\in\hat{\mf{h}}_{red}^*$, $\tilde{M}(\lambda)$ is simple and so $M\cong \tilde{M}(\lambda)$.\\

On the other hand, assume there exists at least one $p_i$ such that $p_i>1$. We are going to construct a set of elements in $M$ which are killed by $e_{i0}$ for all $i\in I_0$. First of all, set $p^{(1)} = \max\{p_i|i\in I_0\}$ and set $w_i:=e_{i0}^{p^{(1)}-1}v$. Note that $w_i=0$ if $p^{(1)}>p_i$ and $w_i\neq 0$ if $p^{(1)} = p_i$, so at least one $w_i$ in non-zero. If for all $j\in I_0$, $e_{j0}w_i=0$ we are done, if not there exists numbers $p_{ij} \in \Z_{>0}$ such that $e_{j0}^{p_{ij}}w_i=0$ and some of the $p_{ij}$ are strictly bigger than $1$. Set $p^{(2)}=\max\{p_{ij} | i,j\in I_0\}$ and set $w_{ij} = e_{j0}^{p^{(2)}-1} w_i$, note that at least one $w_{ij}$ is non-zero. If $e_{k0}w_{ij}=0$ for all $k\in I_0$ we are done, if not we repeat the process. Because of the locally nilpotency of the $e_{l0}$ for $l\in I_0$, in finitely many steps, let say $\ell$ steps, we can find at least one non-zero element $w_{\bf i}$, for ${\bf i} = i_1i_2\ldots i_{\ell}$ a string of elements in $I_0$ such that $e_{l0}w_{\bf i}=0$. Moreover, if ${\bf i}^-$ denotes the string $i_1i_2\ldots i_{\ell -1}$, then $w_{\bf i} = e_{i_{\ell}0}^{p^{(\ell)}-1}w_{{\bf i}^-}$ and so, for all $n\in \Z\setminus \{0\}$, $0=h_{i_{\ell}n}e_{i_{\ell}0}^{p^{(\ell)}}w_{{\bf i}^-} = 2p^{(\ell)}e_{i_{\ell}n}w_{\bf i}$, i.e., $e_{i_{\ell}n}w_{\bf i }=0$. Now, $0 = h_{j0}e_{jm}e_{i_{\ell}0}^{p^{(\ell)}}w_{{\bf i}^-} = e_{jm}h_{j0}e_{i_{\ell}0}^{p^{(\ell)}}w_{{\bf i}^-} + 2e_{jm}e_{i_{\ell}0}^{p^{(\ell)}}w_{{\bf i}^-} = 2p^{(\ell)}e_{jm}e_{i_{\ell}0}^{p^{(\ell)}-1}w_{{\bf i}^-} = 2p^{(\ell)}e_{jm}w_{\bf i}$. \\

Pick one of the non-zero $w_{\bf i}$ constructed above and let $W_{\bf i} = U(G)w_{\bf i}$ be a $G$-submodule of $M$. By construction $e_{ln}W_{\bf i}=0$ for all $l\in I_0$ and $n\in \Z$. Considered the induced module $I(W_{\bf i}) = \Ind_{G\oplus H\oplus N_{+}}^{\hat{\g}} W_{\bf i}$, where $N_{+} = \bigoplus_{i\in I_0, n\in Z}\C e_{in}$ acts by $0$, $H = \bigoplus_{i\in I_0} \C h_i \oplus \C d$ acts by $h_iw_{\bf i} = \mu(h_i)w_{\bf i}$, $dw_{\bf i} = \mu(d)w_{\bf i}$, for some weight $\mu$. Because $M$ is simple, it is a quotient of $I(W_{\bf i})$. If $w_{\bf i}\in T(M)$, we have $W_{\bf i} = \C w_{\bf i}$, and so $M$ is a quotient of $I(W_{\bf i}) = \tilde{M}(\lambda)$ and we are done.\\

In case $w_{\bf i}\notin T(M)$, as in the proof of Proposition 6.0.3. of \cite{CFM03} we get a contradiction. This completes the proof.

\findem

\begin{prop}
    If $M$ is an arbitrary object in $\mc{O}_{red,im}$, then $M\cong \bigoplus_{\lambda_i \in \hat{\mf{h}}^*_{red}} \tilde{M}(\lambda_i)$, for some $\lambda_i's$.
\end{prop}

\dem Because $M$ is in $\mc{O}_{red,im}$, it is a $G$-compatible and so, it has a decomposition as a $G$-module given by $M\cong T(M) \oplus TF(M)$. Since all the weights of $M$ are in $\hat{\mf{h}}^*_{red}$, $T(M)$ is not a $\hat{\g}$-submodule of $M$. Indeed, suppose $T(M)$ is a $\hat{\g}$-module.
let $v\in T(M)$ and consider $f_0 v\in T(M)$. Then $h_{0m}f_0v=0$ and $f_mv=0$ for any $m\neq 0$. Applying $h_{0,-m}$ we get $h_{0,-m}f_m v=0$ and $f_0 v=0$. Since the weight of $v$ is in $\hat{\mf{h}}^*_{red}$, $e_0^p v\neq 0$ for any $p>0$. But if $p$ is sufficiently large the weigh of 
$e_0^p v$ will not be in $\hat{\mf{h}}^*_{red}$ and we get a contradiction. \\

Let $v\in T(M)$ non-zero. As in the proof of the previous statement there exists a string ${\bf i}$ of elements of $I_0$ and a vector $w_{\bf i}$ such that $e_{jm}w_{\bf i} =0$ for all $j\in I_0$ and $m\in \Z$. Let $W_{\bf i} = U(G)w_{\bf i}$. Then we have two possibilities: either $w_{\bf i}\notin T(M)$ or $w_{\bf i} \in T(M)$.\\

In the first case, consider the induced module $I(W_{\bf i})$. Clearly $TF(I(W_{\bf i}))\subseteq I(W_{\bf i})$. Now, if $w\in I(W_{\bf i})$, because $w_{\bf i}\notin T(M)$ we have $gw\neq 0$ for $g\in G$ and so $w\in TF(I(W_{\bf i}))$. Then $TF(I(W_{\bf i})) = I(W_{\bf i})$. By the five lemma, any quotient and subquotient of $I(W_{\bf i})$ also satisfies this property. Set $M' := U(\hat{\g})w_{\bf i}$ which is a subquotient of $I(W_{\bf i})$. Then $M'$ is a $\hat{\g}$-submodule of $M$ and so $M' = TF(M')$ is a $\hat{\g}$-submodule of $TF(M)$, but $TF(M)$ does not have proper $\hat{\g}$-submodule and so $M'=TF(M)$. But, $W_{\bf i}$ is a proper $G$-submodule of $M'$ which is not possible because $M$ is in $\mc{O}_{red,im}$. And so, this case does not occur.\\

In the second case, $W_{\bf i} = \C w_{\bf i} \subseteq T(M)$. So, as $\hat{\g}$-modules $I(W_{\bf i}) \cong \tilde{M}(\lambda_{\bf i})$ for some $\lambda_{\bf i}$, is a $\hat{\g}$-submodule of $M$. Then, any non-zero element of $T(M)$ generates an irreducible reduced imaginary Verma module which is a $\hat{\g}$-submodule of $M$ and because there are no extensions between them, they are direct summands on $M$.

\findem

\begin{cor}
The category $\mc{O}_{red,im}$ is closed under taking subquotients and direct sums, so it is a Serre subcategory.
\end{cor}

\begin{rem} The proofs on the above statements  depends on the structure of reduced imaginary Verma modules, the closed partition $\Delta_{\textrm{nat}}$ and the associated Borel subalgebra ${\mf b}_{\textrm{nat}}$. But, the properties of reduced imaginary Verma modules hold for both untwisted or twisted affine Lie algebras. Moreover, the natural Borel subalgebra for the twisted Lie algebra is properly contained in the natural Borel subalgebra for the untwisted case. So, the results above hold for any affine Lie algebra.
\end{rem}

\section*{Acknowledgement}

JCA has been support by the FAPESP Grant 2021/13022-9. \\

 \bibliographystyle{plain}
 \bibliography{refcatOIVM}

\end{document}